\documentclass[12pt,eqno,sumlimits]{amsart}

\usepackage{amssymb, amscd, amsmath, epsfig, mathtools}
\usepackage[utf8]{inputenc}
\usepackage{amsthm}
\usepackage{enumerate}
\usepackage{xcolor}
\usepackage{scalerel}
\usepackage{soul}
\usepackage{tikz-cd}
\usepackage{esint}
\usepackage{slashed}
\usepackage{hyperref}
\usepackage{mathtools}
\usepackage[new]{old-arrows}

\usepackage[margin=1.0in]{geometry}

\newtheorem{thm}{Theorem}[section]
\newtheorem*{theorem*}{Theorem}
\newtheorem{defn}{Definition}[section]
\newtheorem{cor}{Corollary}[section]
\newtheorem{lem}{Lemma}[section]
\newtheorem{prop}{Proposition}[section]
\theoremstyle{definition}
\newtheorem{rem}{Remark}[section]

\newcommand{\C}{\mathbb C}
\newcommand{\R}{\mathbb R}

\newcommand{\Z}{\mathbb Z}

\newcommand{\calG}{\mathcal G}
\newcommand{\calK}{\mathcal K}

\newcommand{\Psh}{{\rm Psh}}

\newcommand{\diff}{{\rm Diff}}
\newcommand{\ism}{{\rm Isom}}
\newcommand{\conf}{{\rm Conf}}
\newcommand{\cconf}{{\rm conf}}
\newcommand{\psh}{{\rm Psh}}

\newcommand{\est}{\eta, {\rm str}}
\newcommand{\as}{\alpha, {\rm str}}

\newcommand{\sgn}{{\rm sgn}}

\begin{document}

\title[Geometry of Loop Spaces V]{The Geometry of Loop Spaces V: Fundamental groups of
 geometric transformation groups} 
\author[Y. Maeda]{Yoshiaki Maeda}
\address{Tohoku Forum for Creativity, Tohoku University}
\email{yoshimaeda@tohoku.ac.jp}
\author[S. Rosenberg]{Steven Rosenberg}
\address{Department of Mathematics and Statistics, Boston University}
\email{sr@math.bu.edu}
\maketitle

\centerline{In memory of Yuji Ito}  

\begin{abstract}We use differential forms on loop spaces to prove that the fundamental group of certain geometric transformation groups is infinite.  Examples include both finite and infinite dimensional Lie groups. The finite dimensional examples are the conformal group of $S^{4k+1}$ for a family of nonstandard metrics, and the group of pseudo-Hermitian transformations of a compact CR manifold. Infinite dimensional examples include the group of strict contact diffeomorphisms of a regular contact manifold, and other groups coming from symplectic and contact geometry.  
\end{abstract}

\section{Introduction}

In a series of  papers  \cite{MRT3, MRT4, MRT2C, emr, MRTV}, we developed the geometry of loop spaces and a theory of characteristic classes on loop spaces.  In particular, the secondary  Wodzicki-Chern-Simons (WCS) form was used  to prove that $\pi_1(\ism(M,g))$ is infinite for certain closed manifolds $M$ with specific Riemannian metrics $g$.  
In this paper, we 
prove similar results using alternatives to the WCS form, and with the finite dimensional Lie group
$\ism(M,g)$ replaced by certain infinite dimensional Lie groups of geometric transformations.

In our setup, 
$M$ is a closed connected oriented finite dimensional manifold.
Let $\mathcal{G}(M)$ be a closed finite or infinite dimensional subgroup of $\diff(M)$.
We want  examples where $\pi_1(\mathcal{G}(M))$ is infinite.  In general, this seems difficult to prove, even in the explicit case where $\mathcal{G}(M)= \ism(M,g).$
In all our examples, $M$ has an $S^1$ action via diffeomorphisms, and we prove that the associated element of $\pi_1(\mathcal{G}(M))$ has infinite order. 

Our techniques assume a smooth structure on $\diff(M)$ and on the loop space $LM.$  These manifolds have various smooth structures, depending on whether the model space is a Hilbert space \cite{Schmid}, a Banach space \cite{eells}, a Fr\'echet space or a locally convex space \cite{Michor} (of if we take ILH structures \cite{milnoridlg}, \cite{Omori}).  As long as we take the same type of structure on both $\diff(M)$ and $LM$, the results we use from our earlier papers are valid.  If we give $\mathcal{G}(M)$ the induced structure from $\diff(M)$, the results on 
$\pi_1(\mathcal{G}(M))$ are independent of the choice of structure.

Our previous work focused on contact and Sasakian manifolds.
Recall that 
 $(M^{2k+1}, \eta)$ is a contact 
manifold if 
$\eta$ is a one-form on $M$ satisfying
$ \eta \wedge (d\eta)^k \neq 0.$
The characteristic vector field 
(or Reeb vector field) of $(M,\eta)$
is defined by
$d\eta (\xi , \cdot )=0, \eta(\xi) =1.$
A contact manifold $(M, \eta)$
is  regular if the flow of the 
characteristic vector field $\xi$
through any point $m \in M$ is
periodic.  Such manifolds arise as the total space of circle bundles over symplectic manifolds.
The extra conditions which make a contact manifold Sasakian are detailed in \cite{MRTV}.

We proved the following results about the isometry group $\ism(M,g)$ and the strict contactomorphism group 
$\diff_{\est}(M) 
= \{ \phi \in \diff(M), \phi^* \eta = \eta \},$
of contact and Sasakian manifolds $(M,g)$.

\begin{thm}\label{thm:1.1}
(i) \cite[Thm.~3.10]{MRT4},  \cite[Cor.~2.1]{MRT2C} Let $(M,\omega)$ be an integral symplectic manifold of dimension $4k$, and let $\overline{M}_p$ be the total space of the circle bundle with first Chern class $p\omega$.  Then $\overline {M}_p$ admits a Riemannian metric $g_p$ such that for $p\gg 0$, 
$$|\pi_1(\ism(\overline{M}_p, g_p) | = \infty.$$
Equivalently, if $\overline{M}$ is a  regular contact manifold, then $\overline{M}$ covers infinitely many strictly regular contact manifolds $(\overline{M}_p, g_p)$ with 
$|\pi_1(\ism(M_p, g_p) | = \infty.$

(ii) \cite[Thm.~5.1]{MRTV}
Let $(M, g, \phi, \xi, \eta)$ be a connected,
closed, regular $(4k+1)$-dimensional Sasakian manifold.
Deform the Sasakian metric $g$ to the family
$g_\rho= g + \rho^2 \eta \otimes \eta$,
$\rho>0$, where $\eta$ is the
contact $1$-form on $M$. 
Then
\begin{equation*}
 |\pi_1(\ism(M, g_\rho) | = \infty.       
\end{equation*}

(iii) \cite[Thm.~7.1]{MRTV} 
\begin{equation*}
|\pi_1 (\diff_{\est}(S^{4k+1}))| = \infty.
\end{equation*}
\end{thm}
Part (ii) is particularly interesting when $M= S^{4k+1}$ with the standard metric $g= g_0$ and contact structure, since $\pi_1(\ism(S^{4k+1}, g)) = \Z_2$. Therefore, $\pi_1$ of the isometry group is discontinuous in $\rho.$  Note that 
the isometry groups of closed manifolds are finite dimensional Lie groups, 
while the contactomorphism group in (iii) is infinite dimensional.

The main technique in this paper is to replace the WCS $n$-form on the loop space of an  $n$-manifold with other $n$-forms $\hat \calK$ as in (\ref{eq:hatk}).  We give a general condition under which $\hat \calK$ detects an element of infinite order in $\pi_1(\calG(M))$:
\medskip

\noindent {\bf Theorem 2.1}
 {\em Let $M$ be a closed  oriented $n$-dimensional smooth manifold, and let $\calG(M)$ be a subgroup of $\diff(M).$
Let $\hat\calK\in \Lambda^n(LM)$ be defined as in (\ref{eq:hatk}) with kernel $\hat k$ as in (\ref{eq:k}).
If there is a smooth action $a:S^1\times M\to M$  such that (i) $a(\theta,\cdot) \in \calG(M)$ for all 
$\theta\in S^1$, (ii) $a(\theta,\cdot)^*\hat k = \hat k$ for all 
$\theta\in S^1$, (iii) $\int_M a^{L,*}\calK  \neq 0,$
then }
\begin{equation*}
   |\pi_1 (\mathcal{G}(M) )| = \infty.  
\end{equation*}

For example, for one choice of $\hat\calK,$ we obtain a strengthening of 
Thm.~\ref{thm:1.1}(iii).

\medskip
\noindent {\bf Theorem 4.1}
{\em Let $(M, \eta)$ be a closed 
connected regular contact manifold.
Then}  
\begin{equation*}
|\pi_1 (\diff_{\eta, str} (M))| = 
\infty.
\end{equation*}

\noindent This result was previously obtained in \cite{CaSp} using algebraic topology techniques, which gave new information on the cohomology of the classifying space  $B\diff_{\eta, str} (M).$
Our proof is more analytic, and allows us to generalize
Theorem 5.1 
(see \S5.2). Our idea of replacing WCS forms by more general forms on loop spaces is motivated by \cite{CaSp}.

 As an outline of the paper, in \S2 we review some calculations for differential forms on loop spaces, and introduce the forms $\hat\calK$ (\ref{eq:hatk}).  We use these forms  to give a criterion for proving $|\pi_1(\calG(M))| = \infty$ (Thm.~\ref{thm:5.1}). 
 In \S3, we discuss the conformal transformation group
of $S^{4k+1}, k>0.$  We obtain
 
 \noindent {\bf Theorem 3.1} {\em For} $\rho\neq 0$,
 $$|\pi_1(\conf(S^{4k+1}, g_\rho) )| = \infty.$$
 This conformal group is finite dimensional \cite[IV, Thm.~6.1]{kob}.
As above, this theorem fails for the standard metric $g_0.$

 \S4 is devoted to applications of Thm.~\ref{thm:5.1} to finite and infinite dimensional groups of transformations that preserve a contact or canonical one-form. In \S4.1, we discuss the strict contactomorphism group, and a generalization in \S4.2.  In \S4.3, we also prove (Thm.~\ref{thm:5.6})
 $$|\pi_1(\psh(M))| = \infty,$$
  where 
$\psh(M)$ is the finite dimensional group of pseudo-Hermitian transformations of a pseudo-Hermitian (or CR) manifold, the odd dimensional analogue of a symplectic manifold with a compatible almost complex structure.  In \S4.4, we consider a subgroup of the contact transformations $\R^{2k}$, and in \S4.5 we generalize this to the cotangent bundle of a closed manifold. Finally, in \S4.6 we consider Hamiltonian transformations of symplectic manifolds.  
In all cases, we prove that the fundamental group of these transformation groups is infinite by determining the appropriate form on $LM$ to use in Thm.~\ref{thm:5.1}.

We dedicate this paper to the late Professor Yuji Ito 
of Keio University. The first author had the privilege of working alongside Professor Ito as a colleague at Keio University for many years. Despite the
author's lack of prior expertise  in 
ergodic theory and operator algebras, Professor Ito generously and patiently shared his deep knowledge and insights. Professor Ito's intellectual guidance and encouragement
remain a lasting and invaluable asset to the author.

\section{Differential forms on loop spaces and the fundamental group of $\calG(M)$}
In this section, we 
study differential forms on the loop space and some basic properties, as in
\cite{emr, MRTV}. In \S3, we use this material to
give the general method to prove $|\pi_1(\calG(M))|=\infty$ for a
geometric transformation group $\mathcal{G}(M)$. 

Let $M$ be a $n$-dimensional manifold. We consider tensor fields 
${\hat k} \in \Omega^1(M) \otimes 
 \Omega^n (M)$. 
In local coordinates 
$(x^\lambda)= (x^1, \cdots , x^n)$, we have
\begin{equation}\label{eq:k}
{\hat k} = {\hat k}_{\nu [\lambda_1 \cdots , \lambda_n]} dx^\nu \otimes 
dx^{\lambda_1} \wedge \cdots \wedge dx^{\lambda_n}.
\end{equation}
The square brackets denote that the indices are skew-symmetric in
$(\lambda_1, \ldots, \lambda_n)$, and may be omitted if the context is clear.

Let $LM = 
LM =\{ \gamma : S^1 \to M: \gamma \in C^\infty(S^1,M)\}$
be the loop space of $M$.  We fix some $C^\ell$ topology on $LM$ for $\ell\gg 0.$
Given such a   
${\hat k}$, we define
$\hat\calK\in \Omega^n(LM)$ by
\begin{equation}\label{eq:hatk}
{\hat{\mathcal{K}}} (\gamma )
(X_{\gamma ,1}, \cdots , 
X_{\gamma, n}) 
=
\int_{S^1} 
{\hat k}_{\nu [\lambda_1 \cdots \lambda_n]}
(\gamma (\theta)) 
{\dot{{\gamma}}}^\nu (\theta) 
X_{\gamma, 1}^{\lambda_1}(\theta) 
\cdots 
X_{\gamma, n}^{\lambda_n}(\theta) 
d\theta, 
\end{equation}
where 
$X_{\gamma, 1} = X_{\gamma,1}^{\lambda_1}\frac{\partial}{\partial x^{\lambda_1}}, \cdots , X_{\gamma, n}
=X_{\gamma,n}^{\lambda_n}\frac{\partial}{\partial x^{\lambda_n}} \in T_\gamma LM$
are tangent vectors at $\gamma$, 
{\em i.e.,} vector fields along $\gamma \in LM.$  We call $\hat k$ the {\em kernel} of $\hat \calK.$

On an infinite dimensional smooth Banach manifold $N$, the exterior derivative of $\omega\in \Lambda^{s+1}(N)$ is defined by the Cartan formula
\begin{align*}
d_N\omega(X_{0},\cdots X_{s})_p &= \sum_{i=0}^s (-1)^i X_i(\omega(X_{0},\cdots,\widehat{X_i},\cdots, X_{s})\\
&\qquad 
+ \sum_{0\leq i<j\leq s} (-1)^{i+j}\omega([X_i,X_j], X_{0},\cdots, \widehat{X_i},\cdots, \widehat{X_j},\cdots,X_{s}),
\end{align*}
where $X_i\in T_pN$ are extended to vector fields near $p$ using a chart map (see {\it e.g.,} \cite[\S33.12]{KM}). 

We recall the key formulas to show the
infinite order of the fundamental group of geometric transformation group ${\mathcal{G}}(M)$.

As in the finite-dimensional case, we have
 
\begin{lem} \label{lem:A2}
\cite[Lem.~B.2]{emr}
Let $f:W\to L$ be a smooth map between smooth Banach manifolds, and let $\omega\in \Omega^*(L).$
Then $d_W f^*\omega = f^*d_L \omega$. 
 \end{lem}

In fact, the proof carries over to more general settings like ILH manifolds. The following is a consequence of a direct computation
of the exterior derivatives: 

\begin{lem}
\cite[Prop.~B.4]{emr}  
\begin{align}
\label{dLM-formula}
&(d_{LM}{\hat{\mathcal{K}} )(X_{\gamma ,0}, X_{\gamma ,1}, \cdots , X_{\gamma ,2k-1} ) }\\
& \qquad =\sum_{a=0}^n  (-1)^a \int_0^{2\pi}\hat{k}_{\nu \lambda_0 \cdots \hat{\lambda_a} \cdots \lambda_n } (\gamma(\theta))
\dot{X}_{\gamma,0}^{\lambda_0} (\theta) \cdots \hat{X}_{\gamma, a}^{\lambda_a}(\theta)\cdots 
X_{\gamma, n}^{\lambda_n} (\theta) 
d\theta 
\nonumber
\end{align}
\end{lem}

We will consider smooth actions $a:S^1\times M\to M$, with associated maps $a^D:S^1\to\diff(M)$, $a^L:M\to LM$,
given by $a^D(\theta)(x) = a^L(x)(\theta) := a(\theta, x).$  Since we want to study $\pi_1(\calG(M))$ for $\calG(M)\subset \diff(M)$, we need to consider smooth functions $F:[0,1]\times S^1\times M \to M$ and the associated homotopies $F^D:[0,1]\times S^1\to \diff(M),  F^L: [0,1]\times M\to LM$, given by
$F^D(x^0,\theta)(x) = F^L(x^0,x)(\theta) := F(x^0,\theta,x).$

Here we are assuming that $F(x^0,\theta,\cdot)\in \diff(M)$ for all $(x^0,\theta)\in [0,1]\times S^1.$
Then $\{F_*(\partial/\partial x^i)\}_{i=1}^{2k-1}$ is a basis of $T_{F(x^0,\theta,x)}M$ for all $(x^0,\theta,x).$
Therefore, there exist functions $\alpha^i =\alpha^i(x^0,\theta,x)$, $i = 1,\ldots,{\rm dim}(M)$, such that
\begin{equation*}\label{eq:five} F_*\left(\frac{\partial}{\partial x^0}\right) = 
\alpha^iF_*\left(\frac{\partial}{\partial x^i}\right).
\end{equation*}

Using (\ref{dLM-formula}) and replacing of the WCS form $CS^W$  by $\mathcal{K}$ in \cite[Lem.~B.6]{emr}, we have :
\begin{lem}
\label{lem:2.2}
We have
\begin{equation*} 
 F^{L,*} d_{L{\bar M}} \hat \calK (\partial_{x^0}, \partial_{x^1}, \cdots, \partial_{x^{n}}) = 
  \int_0^{2\pi} \hat k_{\lambda_0 [\lambda_1\ldots\lambda_{n}]}
\frac{\partial \alpha^i}{\partial\theta} \frac{\partial F^{\lambda_0}}{\partial x^i}
\frac{\partial F^{\lambda_1}}{\partial x^1}
\cdots
\frac{\partial F^{\lambda_{n}}}{\partial x^{2k-1}}
 d\theta.
 \end{equation*}
\end{lem}

Now we discuss our method for 
proving $|\pi_1(\calG(M))|=\infty,$ for any subgroup $\calG(M)$ of $\diff(M).$  Here $M$ may have nonempty boundary.

We start with some notation.
\begin{defn} (i) Let $f:M\to M$ be smooth. $\hat k$ is $f$-invariant if $f^*\hat k = \hat k.$

(ii) For a smooth map $F:[0,1]\times S^1\times M\to M$, and for
$(x^0,\theta)\in [0,1]\times S^1$,   
set $a(x^0,\theta) = F(x^0,\theta,\cdot): M\to M.$
Then $\hat k$ is $F$-invariant if 
\begin{equation}\label{eq:kinv}
a(x^0,\theta)^* {\hat{k}}   = {\hat{k}},  
\end{equation}    
for all ($x^0,\theta)\in [0,1]\times S^1.$

(iii) $\hat k$ is $\calG(M)$-invariant if $f^*\hat k = \hat k$ for all $f\in \calG(M).$
\end{defn}

For example, the WCS form on a $(4k+1)$-dimensional Riemannian manifold $(M,g)$ is built from
\begin{align} 
\label{K_tensor}
\MoveEqLeft{ \hat k_{\nu [\lambda_1 \cdots \lambda_{4k+1}]}  } \\
&= \sum_{\sigma } \sgn({\sigma})
R_{\lambda_{\sigma (1) e_1 \nu}}{}^{e_{2}} 
R_{\lambda_{\sigma (2) } \lambda_{\sigma (3)} e_{3}}{}^{e_1} 
R_{\lambda_{\sigma (4) } \lambda_{\sigma (5)} e_{1}}{}^{e_{3}} 
\cdots 
R_{\lambda_{\sigma (4k) } \lambda_{\sigma (4k+1)} e_{2}}{}^{e_{k-1}}, \notag
\end{align} 
for $\sigma$ a permutation of $\{1,\ldots,4k+1\}$, and where $R_{ijk}^{\ \ \ \ell}$ are the components of the curvature tensor of  $g$ \cite[App.~B.1]{emr}. 
Then $\hat k$ is $\ism(M,g)$-invariant.

\begin{prop} \label{prop:2.1}  
If $\hat k$ is $F$-invariant, then
\begin{equation*}
 d_{[0,1]\times M}F^{L,*}{\hat{\mathcal{K}}}=0.   
\end{equation*}

\end{prop}

\begin{proof} 
We have 
\begin{equation}\label{eq:5}(a_{(x^0,\theta)}^* {\hat{k}})_{j [i_1, \cdots, i_{n}] }(x) ={\hat{k}}_{\nu [\lambda_1 \cdots \lambda_n]}
(a_{(x^0,\theta)}(x)) 
\frac{\partial a_{(x^0,\theta)}^\nu}{\partial x^j}
\frac{\partial F^{\lambda_1}}{\partial x^{i_1}}
\cdots 
\frac{\partial F^{\lambda_n}}{\partial x^{i_n}},
\end{equation}
where $\partial a_{(x^0,\theta)}^\nu/\partial x^j$ is evaluated at $x$, and the other partial derivatives are evaluated at $a_{(x^0,\theta)}(x) = F(x^0,\theta,x).$
By Lem.~\ref{lem:2.2}, we have
\begin{align*}
F^{L,*} d_{LM} 
\hat{\mathcal{K}}(\partial_{x^0} , \partial_{x^1}, \cdots,\partial_{x^n})
&= \int_0^{2\pi} \frac{\partial \alpha^i}{\partial \theta} 
    \hat{k}_{\lambda_0 [\lambda_1 \cdots \lambda_n]}
\frac{\partial F^{\lambda_0}}{\partial x^i}
\frac{\partial F^{\lambda_1}}{\partial x^1}
\cdots 
\frac{\partial F^{\lambda_n}}{\partial x^n}
d \theta 
\\
&=  \left(\int_0^{2\pi} \frac{\partial \alpha^i}{\partial \theta} 
    d\theta\right) \cdot
\hat{k}_{i [\lambda_1 \cdots \lambda_n]}(x) 
\notag \\
&=0. 
\notag 
\end{align*}
Applying Lem.~\ref{lem:A2} to $F^L:[0,1]\times M\to LM$ gives 
$d_{[0,1]\times M}F^{L,*}{\hat{\mathcal{K}}}=0.$
\end{proof}

We now give a general formulation of \cite[Prop.~3.4]{MRT4}.
\begin{prop} \label{prop:two} 
Let $a_0, a_1:S^1\times M\to M$ be smooth maps such that $a_1^D(\theta), a_2^D(\theta) \in \mathcal{G}(M)$, for all $\theta\in S^1$.

(i)  Let $F:[0,1]\times S^1\times M\to M$ be a smooth homotopy from $a_0$ to $a_1$ with
$a(x^0,\theta)\in \calG(M)$  and  $a(x^0,\theta)(\partial M) \subset \partial M$
for all $(x^0,\theta)\in [0,1]\times S^1$. If $\hat k$ is $F$-invariant,  
then 
$\int_{M} a_0^{L,*}{\hat{\mathcal{K}}}
=  \int_{M} a_1^{L,*}{\hat{\mathcal{K}}}$. 

(ii)  Let 
$a :S^1 \times M \to M$ be a smooth action with $a^D(\theta)\in \mathcal{G}(M)$ for all $\theta\in S^1$.  If 
$\int_M a^{L,*} {\hat{\mathcal{K}}} \neq 0,$ 
then
  $\pi_1(\mathcal{G}(M))$  is infinite.
\end{prop}

\begin{proof} 
(i)  We apply Stokes' Theorem, which is valid for $[0,1]\times M$, which may be a manifold with corners \cite[Thm.~16.25]{lee}.  For $i_{x^0}:M\to [0,1]\times M, i_{x^0}(m) = (x^0,m),$  
we have
\begin{align*}
\int_{M} a_1^{L,*}{\hat{\mathcal{K}}}
-  \int_{M} a_0^{L,*}{\hat{\mathcal{K}}}
&= 
\int_{M} i_1^*F^{L,*}{\hat{\mathcal{K}}}
-  \int_{M} i_0^* F^{L,*}{\hat{\mathcal{K}}}\\
& 
= \int_{[0,1]\times M} d_{[0,1]\times M} F^{L,*}{\hat{\mathcal{K}}}
=0,
\end{align*}
by Prop.~\ref{prop:2.1}.

(ii) Let $a_n$ be the $n^{\rm th}$ iterate of 
 $a$, i.e. $a_n(\theta,m) =
a(n\theta,m).$  
We claim that 
 $\int_{M}a_n^{L,*}{\hat{\mathcal{K}}} =
n\int_{M} a^{L,*}{\hat{\mathcal{K}}}$.  By 
(\ref{eq:hatk}), every term in ${\hat{\mathcal{K}}}$ is of the
form $\int_{0}^{2\pi} \dot{\gamma}(\theta) f(\theta)$, where $f$ is a periodic function on the
circle.  Each loop $\gamma\in
a^L_1(M)$ corresponds to the loop $\gamma(n\cdot)\in a^L_n(M).$  Therefore the term
$\int_{0}^{2\pi} \dot{\gamma}(\theta) f(\theta)$ is replaced by 
$$\int_{0}^{2\pi} 
\frac{d}{d\theta}\gamma(n\theta) f(n\theta)d\theta 
 = n\int_0^{2\pi} \dot{\gamma}(\theta)f(\theta)d\theta.$$
Thus $\int_{M}a_n^{L,*}{\hat{\mathcal{K}}} = n\int_{M}a^{L,*}{\hat{\mathcal{K}}}.$
 By (i), 
$a_n$ and $a_m$ are not homotopic in $\calG(M).$  
By a straightforward modification of \cite[Lem.~3.3]{MRT4}, 
the $[a^L_n]\in \pi_1(\calG(M))$
are all distinct.  
\end{proof}

We now simplify the calculation of $a^{L,*}\hat \calK$ for actions.  For $\hat k\in \Lambda^1(M)\otimes \Lambda^n(M)$ and $\xi\in \Gamma(TM)$, we have the contraction
$\hat k\cdot \xi = \hat k_{j [i_1, \cdots, i_n]}\xi^j\in \Lambda^n(M).$

\begin{lem}\label{lem:2.3}
Let $a: S^1 \times M \to M$ be a smooth action with associated vector field $\xi\in \Gamma(TM)$,
$$\xi^\nu_m = \frac{\partial a^\nu(\theta,m)}{\partial \theta}\biggl|_{\theta = 0}.$$
If $\hat k$ is $a(\theta,\cdot)$-invariant for all $\theta\in S^1,$
then
\begin{equation*}
a^{L,*}\hat \calK = 2\pi \hat k\cdot\xi \in \Lambda^n(M).
\end{equation*}
\end{lem}

\begin{proof}
Since $a$ is an action, we have 
$ a(\theta + \theta^\prime , m) =a(\theta, a(\theta^\prime, m))$, which implies     
\begin{equation*} \xi^\nu_m =
\frac{\partial}{\partial \theta}   
a^\nu (\theta ,m)\biggl|_{\theta=0} = 
\frac{\partial a^\nu}{\partial x^j}\biggl|_{(\theta,m)}
\frac{\partial a^j}{\partial \theta^\prime}{\biggl|_{\theta^\prime =0}}
= \frac{\partial a^\nu}{\partial x^j}\biggl|_{(\theta,m)}\xi^j_m. 
\end{equation*}
Therefore,
\begin{align*}
a^{L,*} \mathcal{K}   
&= \int_0^{2\pi} 
{\hat k}_{\nu [\lambda_1 \cdots \lambda_n]}\xi^j 
\frac{\partial a^\nu}{\partial x^j} 
\frac{\partial a^\lambda_1}{\partial x^{i_1}}
\cdots 
\frac{\partial a^\lambda_n}{\partial x^{i_n}} 
= 2\pi 
\hat k\cdot \xi,
\notag 
\end{align*}
where we write $a(\theta,\cdot)^*\hat k = \hat k$ in local coordinates as in (\ref{eq:5}) to see
that the integrand is independent of $\theta.$

\end{proof}

Combining Prop.~\ref{prop:two}(ii) and Lem.~\ref{lem:2.3} gives the main method to detect
if $|\pi_1(\calG(M))| = \infty.$

\begin{thm}\label{thm:5.1}
 Let $M$ be a closed  oriented $n$-dimensional smooth manifold, and let $\calG(M)$ be a subgroup of $\diff(M).$
Let $\hat\calK\in \Lambda^n(LM)$ be defined as in (\ref{eq:hatk}) with kernel $\hat k$ as in (\ref{eq:k}).
If there is a smooth action $a:S^1\times M\to M$ with associated vector field $\xi$ such that (i) $a(\theta,\cdot) \in \calG(M)$ for all 
$\theta\in S^1$, (ii) $a(\theta,\cdot)^*\hat k = \hat k$ for all 
$\theta\in S^1$, (iii) $\int_M a^{L,*}\calK  = 2\pi \int_M \hat k\cdot\xi
\neq 0,$
then 
\begin{equation*}
   |\pi_1 (\mathcal{G}(M) )| = \infty.  
\end{equation*}    
\end{thm}

\begin{rem}\label{rem:1}
We will use a modified version of this result for regular contact manifolds.  It is easy to check that if we replace (ii) in the Theorem with  $a(\theta,\cdot)^*\hat k = C\cdot  \hat k$ for a nonzero constant $C$, and (iii) with   $\int_M a^{L,*}\calK  = 2\pi C \int_M \hat k\cdot\xi
\neq 0,$  then the proof carries over.
\end{rem}

We will apply this Theorem to various groups $\calG(M)$ in \S\S3-5.  The only real issue is finding a kernel $\hat k$ which is $\calG(M)$-invariant such that (iii) holds.

\section{The Conformal diffeomorphism group of $S^{4k+1}$}

For a Riemannian manfold $(M,g)$,  
the group of the conformal diffeomorphisms of $(M,g)$ is
\begin{equation*}\label{conf2}
 \conf(M, g) =\{ \psi \in \diff(M) | 
 \psi^* g = f g \
 \text{for some}\ f \in C^\infty (M), f>0 \}. 
\end{equation*}

Let $g_{st}$ be the standard metric on $S^{4k+1}$. For the Hopf fibration $\pi: S^{4k+1} \to \mathbb{CP}_{2k}$, the unit vector field $\xi$ along the fiber is  
the Reeb vector field for the standard contact structure. 
The contact one-form $\eta$  
is the dual of $\xi$:
\begin{equation*}
 \eta (\xi) =1,  \  d\eta (\xi , \cdot )=0.    
\end{equation*}
For a real parameter  $\rho \geq 0$, 
we take a new metric on $S^{4k+1}$
\begin{equation*}  
 g_\rho :=  g_{st} +\rho^2 \eta \otimes \eta.
\end{equation*}

In \cite[Thm.~6.1]{MRTV}, we proved that  the $|\pi_1(\ism(S^{4k+1},g_\rho))| = \infty $ iff $\rho>0.$
In this section, we study $\pi_1(\conf(S^{4k+1}, g_\rho))$
by choosing an appropriate kernel in (\ref{eq:k}).  
Namely, we have 
\begin{thm}\label{thm:conf}
\begin{equation*}
| \pi_1 (\conf(S^{4k+1},g_\rho))| = \infty,      
\end{equation*}
for $\rho \neq 0.$
\end{thm}

This results fails if $\rho =0$: $\conf(S^{n-1}, g_{st})$ is diffeomorphic to $SO(n-1,1)$, which has the homotopy type of its maximal compact subgroup $SO(n-1).$  Thus $\pi_1(\conf(S^{n-1}, g_{st})) \simeq \Z_2,$ a type of ``discontinuity" in the fundamental group as $\rho\to 0.$

\begin{proof} We verify the three conditions in Thm.~\ref{thm:5.1} for $\calG(S^{4k+1}) = \conf(S^{4k+1},g_\rho)$.  For the action $a$, we take $a(\theta,\cdot)$ to be the rotation by angle $\theta$ in the circle fibers of $\pi.$  This is an action by isometries \cite[Cor.~4.1]{MRTV}, so it is an action by conformal diffeomorphisms.  Thus (i) holds.

To define $\hat k$ and verify (ii), we compute the Weyl tensor of $g_\rho.$
Recall that $S^{4k+1}$ has the standard Sasakian structure $(g_{st}, \phi, \xi,\eta)$ where $\phi$ is the odd dimensional analogue of an almost complex structure: $\phi_i{}^k \phi_k{}^j
        = -\delta_i{}^j + \eta_i \xi^j$ \cite[\S2]{MRTV}.
Let $R_{kji}{}^h$ and $\bar{R}_{kji}{}^h$ be 
the curvature tensors of $g=g_{st}$ and $g_{\rho},$ 
 respectively.
By \cite[Lem.~4.3]{MRTV},
\begin{align}\label{eq:a}
R_{kji}{}^h   &= g_{ki} \delta_j{}^h - g_{ji}\delta_k{}^h  \\
\bar{R}_{kji}{}^h &= R_{kji}{}^h -\rho^2( \phi_{ki}\phi_{j}{}^{h} -\phi_{k}{}^{h} \phi_{ji}
+2\phi_{kj} \phi_{i}{}^{h} +2\eta_k \eta_i \delta_j{}^h -2 \eta_j \eta_i \delta_k{}^{h} 
+g_{ki}\eta_j \xi^h -g_{ji}\eta_k \xi^h)\notag
 \\
&\qquad -\rho^4(\eta_k \eta_i \delta_j{}^{h} 
-\eta_j \eta_i \delta_k{}^{h} ). 
\notag
\end{align}
This implies
\begin{align}\label{eq:b}
\bar{R}_{ji} &=R g_{ji}-\rho^2(\phi_{ki}\phi_{j}{}^{k} -\phi_{k}{}^{h} \phi_{ji}+2\phi_{kj} \phi_{i}{}^{k} +2\eta_i \eta_j -2(4k+1) \eta_i \eta_j+\eta_j \eta_i -g_{ji})\notag \\
&\qquad -\rho^4(\eta_j \eta_i -(4k+1) \eta_j \eta_i) \\
&= (R^2+4\rho^2) g_{ij} +((4k-1) \rho^2 +4k\rho^4)\eta_j \eta_i,\notag
\end{align}
where $\bar{R}_{ji}= \bar{R}_{kji}{}^k, \bar{R} = h^{ji}\bar{R}_{ji}$ are the Ricci tensor and scalar curvature of $g_\rho$, respectively, and $R = (4k+1)(4k)$ is the scalar curvature of $g$.
The Weyl curvature tensor for $\bar{R}$ is
\begin{align}\label{eq:Weyl}  \bar{C}_{kji}{}^h 
&= \bar R_{kji}{}^h +\frac{1}{4k-1} (\bar{R}_{ki} \delta_j{}^h  -\bar{R}_{ji}\delta_k^h 
  +g_{ki} \bar{R}_j{}^h -g_{ji}\bar{R}_k{}^h)  \\
&\qquad - \frac{\bar{R}}{(4k)(4k-1)}     (g_{ki} \delta_j^h -g_{ji}\delta_k{}^h ). \notag
\end{align}
Plugging 
(\ref{eq:a}) and (\ref{eq:b})
into (\ref{eq:Weyl}, we have
\begin{align}\label{eq:ten}
\bar{C}_{kji}{}^h 
&= -\rho^2
(\phi_{ki}\phi_{j}{}^{h} 
-\phi_{k}{}^{h} \phi_{ji}
+2\phi_{kj} \phi_{i}{}^{h} )
\\
&+c_1(g_{ki} \delta_j{}^h - g_{ji}\delta_k{}^h)
+ c_2 (\eta_k \eta_i \delta_j{}^h 
- \eta_j \eta_i \delta_k{}^{h} )
+c_3 (g_{ki}\eta_j \xi^h -g_{ji}\eta_k \xi^h) 
\notag 
\end{align}
where $c_1, c_2, c_3$ are explicit nonzero constants depending
on $\rho$. 

We set 
\begin{equation*}\hat{k}^{\cconf}_{j i_1 \cdots i_{4k+1}} 
=\bar{C}_{i_1 \ell_1 j}{}^{\ell_0}
\bar{C}_{i_2 i_3 \ell_2}{}^{\ell_1}
\cdots 
\bar{C}_{i_{4k} i_{4k+1} \ell_0}{}^{\ell_{2k}}.
\end{equation*}
This is the conformal version of (\ref{K_tensor}); the similar expressions \cite[(17), (18)]{MRTV}  were used to prove
$|\pi_1(\ism(S^{4k+1}, g_\rho))| = \infty$ in \cite[Thm.~6.1]{MRTV}.
Since the Weyl tensor is invariant under conformal transformations of the metric, so is $\hat k^{\cconf}$.
Therefore, (ii) is verified:
\begin{equation*}  a(\theta,\cdot)^* \hat k^{\cconf} = \hat k^{\cconf}.
\end{equation*}

On the loop space $LM = \{ \gamma: S^1 \longrightarrow M\}$, we consider tangent vectors 
$X_{1}, \cdots, X_{4k+1} \in
 T_{\gamma}(LM) = \Gamma(\gamma^*TM)$ and define $\calK^\cconf\in \Lambda^n(LM)$ by 
\begin{equation}
{\mathcal{K}^{\cconf}}
(X_1,\ldots, X_n)_\gamma
=\int_0^{2\pi} 
\hat{k}^{\cconf}_{j[i_1 \cdots i_{4k+1}]}(\gamma (\theta) )\dot\gamma^j(\theta)
 X_{1}^{i_1}(\theta) \cdots X_{ 4k+1}^{i_{4k+1}}   \ d\theta 
\end{equation}
In \cite[Prop.~6.1]{MRTV}, we computed that for $M = S^{4k+1},$ we have $\hat k = C_\rho \eta\wedge (d\eta)^{2k}$ for 
$\hat k$ in (\ref{K_tensor}) and for some nonzero constant $C_\rho$. The Weyl tensor has the same symmetries as the Riemann curvature tensor, and the terms in (\ref{eq:ten}) are the same as in
\cite[(24)]{MRTV}, so the proof of \cite[Prop.~6.2]{MRTV} carries over to $\hat k^\cconf.$  Thus 
$\int_M a^{L,*}\calK^\cconf \neq 0$ by Lem.~\ref{lem:2.3}.  This verifies (iii).
\end{proof}

\section{Applications to Geometric transformation groups preserving one-forms}

In this section, we discuss the 
the fundamental groups of  geometric transformation groups which preserve certain one-forms.  In
\S4.1,  we prove that the group of strict contactomorphisms has infinite fundamental group (Thm.~\ref{thm:4.1}). In \S4.2, we  relax the conditions on the contact one-form to prove similar results for other groups of diffeomorphisms. In \S4.3, we discuss pseudo-Hermitian transformations.  In \S4.4, we consider transformations of $\R^{2k}$ which preserve a standard one-form, and in \S4.5 we generalize this to the cotangent bundle of a closed manifolds. The groups in these subsections are infinite dimensional, except in \S4.3.

\subsection{The group of strict contact transformations}

Let $(M, \eta)$ be a $(2k+1)$-dimensional connected closed contact manifold, where
$\eta$ is the contact one-form.  We assume that 
$(M, \eta)$ is regular, {\em i.e.}, its  
Reeb vector field $\xi$, characterized by 
\begin{equation}\label{xi}
 \eta(\xi) =1, \ d\eta (\xi, \cdot )=0,  
\end{equation}
has closed orbits.

Let 
\begin{equation*}
\diff_{\est} (M) 
=\{ \phi \in \diff(M); \phi^* \eta = \eta \}
\end{equation*}
be the group of strict contactomorphisms.

\begin{thm}\label{thm:4.1}
Let $(M,\eta)$ be a $(2k+1)$-dimensional closed regular contact manifold. 
Then 
\begin{equation}
|\pi_1(\diff_{\est} (M))| = \infty.   
\end{equation}
In particular, the homotopy clsss of 
the Reeb flow in 
$\pi_1(\diff_{\est}(M))$ has  
infinite order. 
\end{thm}

As mentioned in the Introduction, this was first proved in \cite{CaSp}.

\begin{proof}

As in \S2, we  set 
\begin{align*} 
\hat{k}^{\eta}   
&= 
\eta\otimes \left(\eta \wedge (d\eta)^k\right),   \\
\hat{\mathcal{K}}^{\eta}
(X_{\gamma, 1}, \cdots , X_{\gamma, 2k+1})_\gamma
&= \int_0^{2\pi} 
\hat{k}^{\eta}(\gamma(\theta))_{j[i_1 \cdots i_{2k+1}]}
X_{ 1}{}^{i_1} (\theta) 
\cdots X_{2k+1}{}^{i_{2k+1}}(\theta)d\theta \in \Lambda^{2k+1}(LM).
\end{align*}
where $\gamma \in LM$ and 
$X_{ 1}, \cdots, X_{ 2k+1} \in T_{\gamma} (LM)$.

We now find an action $a_{\rm rot}:S^1\times M\to M$ satisfying Conditions (i) -- (iii)  in Thm.~\ref{thm:5.1}, as modified in Rem.~\ref{rem:1}.  This will complete the proof.

Let $\psi_t (m)$ be the one-parameter group generated by $\xi.$ 
After choosing a metric associated with the contact structure (see \cite[\S2]{MRTV}),  
we can apply Wadsley's theorem \cite{Wadsley} to conclude that there exists $N>0$ such that $N$ is an integral multiple of the period of each Reeb orbit. Therefore, we can modify the flow to $\bar\psi_t(m) := \psi_{(2\pi)^{-1} N t}(m)$ to get an $S^1$ action  
\begin{equation*}
a_{\rm rot}: S^1 \times M \to M,\ a_{\rm rot}(\theta , m) = \bar\psi_\theta (m).
\end{equation*} 
It follows from (\ref{xi}) and the Cartan formula for the Lie derivative that
 $\mathcal{L}_{\xi} \eta =0. $  
As in Rmk.~\ref{rem:1},   
$a_{\rm rot}(\theta,\cdot)^* \eta =(N/2\pi)\eta$ (Condition (i))
which implies $a_{\rm rot}(\theta,\cdot)^*\hat k^\eta = C\cdot\hat k^\eta$ for some $C\neq 0.$ (Condition (ii)).
For Condition (iii), by Lem.~\ref{lem:2.3},
\begin{equation*}
\int_M (a_{\rm rot}^L)^* {\hat{\mathcal{K}}} 
= 2\pi C \int_M \hat k^\eta\cdot\xi = 2\pi C\int_M \eta(\xi)\ \eta \wedge (d\eta)^k
=2\pi C\int_M \eta \wedge (d\eta)^k  \neq 0.
\end{equation*}
\end{proof}

\subsection{Generalizations of contactomorphism groups}

The proof of Thm.\ref{thm:4.1} immediately carries over to more general setups.
Let $M$ be a closed, connected, oriented smooth $(2k+1)$-manifold, and let
$\eta$ be a one-form on $M$. Assume  there is a vector field $\xi$
on $M$ that satisfies the following:
\begin{enumerate}
\label{assumption1.3}
\item [(A1)]
The flow of the vector field $\xi$ is periodic with period independent of the orbit. 
\item [(A2)]
$L_{\xi} \eta =0. $
\item [(A3)] 
$\int_M \eta(\xi) \eta \wedge (d \eta)^{2k} \neq 0.$
\end{enumerate}
Set 
$\diff_{\eta} (M) = \{\phi\in \diff(M): \phi^*\eta = \eta\}.$

\begin{thm}
 Under the assumptions (A1)-(A3), we have
 \begin{equation*}
|\pi_1 (\diff_{\eta} (M))| = 
\infty.
\end{equation*}
\end{thm}

We  give an example satisfying 
(A1) -- (A3).
Let $T^3 = S^1 \times S^1 \times S^1$ be the 3-torus with 
 coefficients $u = (u^1, u^2, u^3)$. 
Set
\begin{equation*}
 \eta (u) = \eta_1 (u^2) du^1 + \eta_3 (u^2) du^3, 
 \  \xi (u) = \partial_{u^1} .
\end{equation*}
Clearly, (A1) holds, and it is easily checked that $L_{\xi}  \eta =0.$
Noting that
\begin{equation*}
d\eta = \frac{\partial \eta_1}{\partial u^2} du^2 \wedge du^1
+ \frac{\partial \eta_3}{\partial u^2} du^2 \wedge du^3,    
\end{equation*}
we have
\begin{equation*} \int_{T^3} \eta (\xi) \eta \wedge d\eta = (2\pi)^2 \int_0^{2\pi} 
\eta_1(u^2)\left(\eta_1(u^2) \frac{\partial \eta_3}{\partial u^2} -
\eta_3 (u^2)\frac{\partial \eta_1}{\partial u^2} \right)du^2.
\end{equation*}
For $\eta_1 (u^2) = \cos u^2 , \eta_3 (u^2) = \sin 2 u^2,$ we get 
\begin{align*}
\int_{T^3} \eta (\xi) \eta \wedge d\eta &= (2\pi)^2 \int_0^{2\pi} \left(2 (\cos u^2)^2 \cos 2u^2
-\cos u^2 \sin 2 u^2 \sin u^2\right)du^2= 2\pi^3\neq 0.
\end{align*}

Thus, we get
\begin{cor}
Let $T^3$  be a 3-torus with coordinates $(u^1, u^2,u^3)$,
let 
$\eta (u) = \eta_1(u^2) du^1 + \eta_3 (u^2) du^3$, 
and let $\xi = \partial_{u^1}.$
Then
$|\pi_1 (\diff_{\eta} (T^3))|= \infty$.  
Specifically, the loop of diffeomorphisms given by rotation in the $u^1$ direction has infinite order in $\pi_1 (\diff_{\eta} (T^3))$.
\end{cor}

In a second direction, we can replace $\eta\wedge (d\eta)^k$ with a general top degree form.
Let $M$ be an oriented closed $C^\infty$ $n$-manifold.  We choose 
$\eta\in  \Lambda^1 (M)$ and $\mu \in \Lambda^n (M)$.  
We note that $\mu$ is not necessarily a volume form on $M$.

Set $\diff_{\mu, \eta }(M)  =\{\phi \in \diff(M): \phi^*\eta = 
\eta, \phi^*\mu = \mu.\}$.

  If $\mu$ is a volume form, then the Lie algebra of $\diff_{\mu} (M) = \{\phi\in\diff(M): \phi^*\mu = \mu\}$ is the space of divergence-free vector fields, which is infinite dimensional. We expect that $\diff_{\eta} (M)$ and $\diff_{\mu, \eta } (M) $ are also infinite dimensional.

We assume that there is a vector field $\xi$ on $M$ such 
\begin{enumerate}
\label{assumption1.3a}
\item [(B1)]
The flow of the vector field $\xi$ is periodic with period independent of the orbit.  
\item [(B2)]
$L_{\xi} \eta =  L_\xi \mu = 0. $
\item [(B3)] $\int_M \eta (\xi ) \mu \neq 0$.
\end{enumerate}

Then, we have 
\begin{thm}
\label{mu-theorem}
 Under the  assumptions (B1) -- (B3), we have 
$|\pi_1 (\diff_{\mu, \eta} (M)) | = \infty$. 
\end{thm}
\begin{proof}
We take 
\begin{equation*}
\hat k = \eta\otimes \mu \in 
\Lambda^1 (M) \otimes \Lambda^n (M).    
\end{equation*}
It is clear that $\hat k$ is invariant under the group
$\diff_{\mu, \eta} (M)$. 
Since $ \int_M \hat k\cdot \xi = \int_M {\hat{k}}_{j [i_1 \cdots i_n] } \xi^j \neq 0$, Lem.~\ref{lem:2.3} and Thm.~\ref{thm:5.1} give the result.
\end{proof}

We give a simple example that Theorem \ref{mu-theorem} holds on the torus
$T^2 = S^1 \times S^1$ with coordinates $(u^1, u^2)$. 
The flow of $\xi = \partial_{u^1} $ satisfies (B1).
Set
\begin{equation}\label{eq:mueta}
\mu = du^1 \wedge du^2, \ \eta(u) = \eta_1 (u^2) du^1 + \eta_2 (u^2) du^2,    
\end{equation}
where $\eta_1(u^2) >0$ on $T^2$. 
Then 
$L_\xi \eta =0$, and $L_\xi \mu =0$. 
Note that 
\begin{equation*}
 \int_{T^2} \eta (\xi ) \mu = 
 \int_{T^2} \eta_1 (u^2) du^1 \wedge du^2  
 =2\pi \int_0^{2\pi} \eta_1(u^2) du^2 > 0.
\end{equation*}
Thus, (B1) -- (B3) are satisfied, and we have
\begin{cor}
For the choice of $\mu$ and $\eta$ in (\ref{eq:mueta}),  we have
$|\pi_1 (\diff_{\mu, \eta} (T^2))| = \infty$. 
\end{cor}

\subsection{The group of pseudo-Hermitian transformations}
We discuss 
the transformation group of a psuedo-Hermitian structure
(CR structure) on a closed regular contact manifold. 
Let $M$ be a closed regular contact manifold with 
contact form $\eta$. Assume that there exists a complex
structure $J$ on the contact bundle $\text{Ker}\  \eta$ and that the
Levi form $d\eta \circ J$ is a positive definite 
Hermitian form.  Then, $( \eta, J)$ is called a 
psuedo-Hermitian structure on $M$. For the
Riemannian metric 
$    g := d\eta \circ J + \eta \otimes \eta$
on $M$,  the group of pseudo-Hermitian 
transformations of $M$ is 
\begin{equation*}
\Psh(M) = \{ h \in \diff(M) \, | \, h^* \eta = \eta, 
h_* \circ J = J \circ h_*: \text{ on } \text{Ker}\,  \eta  \}.
\end{equation*}
We note that 
\begin{equation}
\label{psh_isom}
  \Psh(M) \subset \ism(M, g).
\end{equation}

We have the following result. 
\begin{thm}\label{thm:5.6}
Let $M$ be a $(2k+1)$-dimentional closed regular contact
manifold with a psuedo-Hermitian structure. 
Assume that the Reeb vector field $\xi$ defined by
$\eta (\xi) =1$ and $d\eta (\xi , \cdot)=0$ generates
a periodic one-parameter transformation group 
of psuedo-Hermitian transformations.
Then $|\pi_1(\Psh(M)|= \infty$. 
\end{thm}
\begin{proof}
We take as kernel function 
\begin{equation}
\label{kernel_psh}
\hat k^{\Psh} = \eta\otimes {\rm dvol}_g,
\end{equation}
where ${\rm dvol}_g$ is the volume form for $g$. 
By (\ref{psh_isom}), the kernel function
defined by (\ref{kernel_psh}) is preserved 
by  pseudo-Hermitian transformations. 
Note that
\begin{equation*}
{\mathcal{K}}^{psh} = \int_{M}  
    \hat k^{\Psh}\cdot\xi = {\rm vol}(M) 
   \neq 0.
\end{equation*}
Thus, 
 Thm.~\ref{thm:4.1} gives the result.
\end{proof}

\subsection{The group of symplectic transformations of homogeneous degree one on $\mathbb{R}^{2k}$ }
Let $\mathbb{R}^{2k}$ be Euclidean
$2k$-space with the one-form
\begin{equation*}
\alpha = \frac{1}{2} 
  \sum_{i=1}^k x^i d\xi^i -\xi^i dx^i,    
\end{equation*}
where $z=(z^1, \cdots, z^k)$ are the coordinates on 
$\mathbb{R}^{2k}= \C^k$ and $z^i=(x^i, \xi^i)$. We note that $d\alpha = \omega$ is the standard symplectic form on $\R^{2k}.$

For $\overset{\circ}{\mathbb{R}}{}^{2k}= \mathbb{R}^{2k} - \{0\}$,
$\phi \in \diff(\mathbb{R}^{2k})$ is of homogeneous degree one if 
\begin{equation*}
 \phi (rx,r\xi) = r\cdot \phi(x,\xi) \,  
  \text{ for } r>0.    
\end{equation*}

Let $\diff^{(1)} (\mathbb{R}^{2k})$ be
the subgroup of 
$\diff(\mathbb{R}^{2k})$ consisting of homogeneous degree one diffeomorphisms.
$\diff^{(1)}(\mathbb{R}^{2k})$ contains the 
subgroup
\begin{align*}
\diff^{(1)}_{\as} (\mathbb{R}^{2k}) 
&= \{ \phi \in \diff^{(1)} (\mathbb{R}^{2k}) \,| 
  \phi^* \alpha = \alpha\}.
\end{align*}
This group is an infinite dimensional
Lie group with good differential structures
e.g., ILH-structures, Fr\'echet structures, etc. \cite{Omori}.

Let $S^{2k-1}$ be the unit sphere of $(2k-1)$ dimensional with the origin as the center, and
$i: S^{2k-1} \longrightarrow \mathbb{R}^{2k}$ be the standard 
embedding of the unit sphere $S^{2k-1}$ into $\mathbb{R}^{2k}$.
We define a map $A:\diff(S^{2k-1})\to \diff^{(1)}(\mathbb{R}^{2k})$ by
$(A\hat\phi)(r\hat x, r\hat\xi) = r \hat\phi(\hat x, \hat \xi)$, for $(\hat x, \hat \xi)\in S^{2k-1}$
and $\hat\phi\in \diff(S^{2k-1}).$ $A$ is clearly not surjective. We define
\begin{align*}
\diff^{(1)}_{\as, A} (\mathbb{R}^{2k})
&=
\{ \phi \in \diff^{(1)}_{\as}(\mathbb{R}^{2k}) \, | \,
\phi\in {\rm Im(A)} \}.
\end{align*}

Note that $\bar{\alpha}: = i^* \alpha$ 
gives a contact structure on $S^{2k-1}$. 
Let $\phi \in \diff^{(1)}_{\as,A}(\mathbb{R}^{2k})$.
We have
\begin{equation*}
\phi^*(r\alpha) = r\hat{\phi}^*(\bar{\alpha}) ,    
\end{equation*}
where
$\hat{\phi}(\hat{x}, \hat{\xi})= \phi \circ i(\hat{x}, \hat{\xi})$. Thus,
$\hat{\phi}  \in \diff_{\bar{\alpha},\text{str}}(S^{2k-1})$, which implies

\begin{lem}\label{lem:5.1}
\begin{equation*}
\diff^{(1)}_{\as, A}(\mathbb{R}^{2k}) =
\diff_{\bar{\alpha}, {\rm str}} (S^{2k-1} ).
\end{equation*}
\end{lem}

Since $S^{2k-1}$ is a regular contact manifold, Thm.~\ref{thm:4.1} implies
\begin{cor}\label{cor:5.1}
\begin{equation*}
|\pi_1(\diff^{(1)}_{\as, A} (\mathbb{R}^{2k} ))| 
= \infty.
\end{equation*}
\end{cor}

\subsection{The group of canonical transformations of degree one
on the cotangent bundle}

Considering $\R^{2k}$ as $T^*\R^k$, we have a similar situation for the cotangent bundle of a symplectic manifold.
Let $(M^k,\omega) $ be a closed, connected $C^\infty$ symplectic manifold with cotangent bundle $\pi:T^*M \longrightarrow M.$
On  $\overset{\circ}{T^*M} := T^*M - \{\text{zero-section}\}$, we define the canonical/contact one-form $\alpha = \sum_{i=1}^k \xi^i dx^i $
in local coordinates $(x, \xi)$ on $T^*M$.

Let 
$\diff^{(1)} (\overset{\circ}{T^*M})$ be the group of
diffeomorphisms of $\overset{\circ}{T^*M}$ of homogeneous
degree 
one in  
the fiber direction; {\em i.e.,}
if we write $\phi \in \diff^{(1)} (\overset{\circ}{T^*M})$
as
\begin{equation*}
\phi (x,\xi) = (\phi_1 (x,\xi ) , \phi_2 (x,\xi)),
\end{equation*}
where $\phi_1,$ resp.~$ \phi_2$, involve only $x $, resp.~$\xi$, coordinates,
then
\begin{equation}\label{phi22}
 \phi_2(x, r\xi) = r\phi_2(x, \xi) \text{ for } r>0.    
\end{equation}
Since $\phi_2$ changes by a function of $M$ only under a change of coordinates on $M$, 
(\ref{phi22}) is independent of local coordinates.

We set 
\begin{equation*}
\diff^{(1)}_{\as} (\overset{\circ}{T^*M}) 
= \{ \phi \in \diff^{(1)} (\overset{\circ}{T^*M})
\, | \phi^* \alpha =\alpha   \}. 
\end{equation*}
This 
is also an infinite dimensional Lie group, since 
for a $C^\infty$ diffeomorphism 
$f: M \longrightarrow M$, we  
have (\ref{phi22}) for $(df)^*$ (the adjoint of the differential $df$), and $(df)^{*} \alpha = \alpha$  by the cotangent bunlde lift theorem \cite[Prop.~6.3.2]{Marsden-Ratiu}. (Here we abuse notation by using $(df)^*$ instead of $(df)^{**}$ for the pullback on one-forms associated to $(df)^*$.)

A choice of metric $g$ on  $M$ gives an inner product on each cotangent fiber and allow us to define
the unit cosphere bundle
$S^*M= \{ (x, \xi) \in \overset{\circ}{T^*M}| \ |\xi|_g =1 \}$.  We note that for the inclusion $i:S^*M\longrightarrow T^*M$, $\bar\alpha := i^*\alpha$ is a contact form on $S^*M$.
We set
\begin{equation*}
\diff^{(1)}_{\as, g} (\overset{\circ}{T^*M}) 
= \{ \phi \in \diff^{(1)}_{\as} (\overset{\circ}{T^*M})
| \  |(x,\xi)|_g = 1\Rightarrow| \phi_2(x,\xi) |_g=1 \}.    
\end{equation*}

Let $a:S^1\times M\to M$ be a smooth $S^1$-isometric action: 
 {\em i.e.,} $a$ is a smooth action, and $a^D:S^1 \to \diff(M)$, defined by $ a^D(\theta)(x) :=
a(\theta,x),$ has $a^D(\theta)\in \ism(M,g)$ for all $\theta$. Then we have:

\begin{cor}
Let $(M^k,g,\omega)$ be a closed connected $C^\infty$ Riemannian  $k$-manifold with 
a smooth $S^1$-isometric 
action 
on $(M,g)$. Then 
\begin{equation*}  
|\pi_1 (\diff^{(1)}_{\as,g} (\overset{\circ}{T^*M})|
= \infty.   
\end{equation*}
\end{cor}
\begin{proof}
 Since $a^D$ is an $S^1$ action on $\diff(M),$ the adjoint $(da^D)^*$ of the differential $da^D$ is an $S^1$ action on $T^*M$.  Since $(da^D)^*$ is linear in each fiber, $(da^D)^*$ gives an action on 
 $\diff^{(1)}_{\alpha,A} (\overset{\circ}{T^*M})$; this uses
 $(da^D)^*\alpha = \alpha$ as above.  
 As in Lem.~\ref{lem:5.1}, 
 $$\diff^{(1)}_{\as,g} (\overset{\circ}{T^*M}) = \diff_{\bar \alpha, {\rm str}}(S^*M),$$
 where $S^*M$ is the unit cotangent bundle and $\bar\alpha = i^*\alpha$ for the inclusion $i:S^*M \longrightarrow T^*M.$
 
For
 each $\theta$, 
$|(x,\xi)|_g = 1$ implies $|da^D(\theta)(x.\xi)|_g = 1$, since the action is via isometries. Therefore, $(da^D)^*$ descends to an $S^1$ action  $\hat a$ on
$\diff_{\bar \alpha, {\rm str}}(S^*M)$ satisfying 
condition (iii) in Thm.~\ref{thm:5.1}:
$\int_{S~*M} \hat{a}^{L,*}\mathcal{K} \neq 0,$
where 
 $\mathcal{K}$ has kernel
$k= \bar{\alpha} \otimes \bar{\alpha} \wedge (d\bar{\alpha})^k. $
By Thm.~\ref{thm:5.1}, $\pi_1(\diff_{\bar \alpha, {\rm str}}(S^*M))$ is infinite, which gives the result.
\end{proof}

\subsection{Lie group of Hamiltonian symplectic transformations} 

We give an application to  
an interesting subgroup of the Lie algebra of the Poisson algebra of smooth functions on a 
symplectic manifold $N$. The main reference for this subsection is \cite{Visman}.

Let $(N,\omega)$ be a closed symplectic $2k$-dimensional manifold 
with a symplectic form $\omega$, which is 
is integrable, {\em i.e.,} 
$[\omega] \in H^2(N, \mathbb{Z}).$ 
Then there is
an $S^1$-bundle $\pi: (M, \eta) \longrightarrow (N,\omega) $, 
where $(M, \eta)$ is a contact manifold  
with $\pi^* (\omega) = d\eta$. 
As usual, for a smooth
function $H(x,\xi)$ on $(N,\omega)$, we define the Hamiltonian 
vector field $X_H$ by 
$ \omega( X_H, \cdot)   = dH,$
and define the Poisson bracket $\{ \ , \ \}$ by
$  \{ H , H^\prime \} = X_H H^\prime.$  
It is standard that $(C^\infty (N), \{ \, , \, \})$ is
an infinite-dimensional Lie algebra. 

We consider a  
vector field $ V_f$ on $(M, \eta)$ associated with a smooth function $f \in C^\infty (M)$ defined by
\begin{equation*}
 \eta (V_f) = -f , \quad 
d\eta ( V_f , \cdot )  = df 
\end{equation*}
It is easily seen that $L_{V_f} \eta = 0$, so $V_f$ is by definition a strict contact vector field. 

Let $\mathcal{X}_{\eta , {\rm str}}(M) $ be the Lie algebra
of  strict contact vector fields on 
$(M, \eta)$.
For any $V \in \mathcal{X}_{\eta , st} (M)$,
there is a smooth function $f$ such that
$V=V_f$, so 
\begin{equation*}
\mathcal{X}_{\eta, {\rm str}} (M) = \{ V_f | f \in C^\infty (M) 
\}. 
\end{equation*}
For $H \in C^\infty (N)$, we denote by  
$H^L \in C^\infty (M) $ the lift
of $H$, {\em i.e.,} $H^L = \pi^* H$,
and set
\begin{equation*}
\mathcal{X}_0 (M) 
  = \{ V_{H^L} | H \in C^\infty (N) \}.    
\end{equation*}
It is easily seen that $\mathcal{X}_0 (M)$
is a closed Lie algebra of 
$\mathcal{X}_{\eta ,{\rm str}}.$

We define 
the contact diffeomorphism
$\phi = \phi_H = \exp( V_{H^L})$ on $(M, \eta)$.
We set $G_0(M)$ to be the Lie group which is  finitely generated
by $\exp (V_{H^L})$, and let $\bar{G}_o(M)$ be the
closure of $G_0(M)$ in $\diff_{\alpha , {\rm str}} (M)$.

This procedure gives a  Lie group $\bar{G}_0(M)$, which we call the Lie group of Hamiltonian symplectic transformations (cf.~\cite{Visman}), whose Lie algebra is a subalgebra   of $(C^\infty (N), \{ \cdot , \cdot \})$.

\begin{thm}
Let $(N, \omega)$ be an integral closed symplectic manifold. Then
\begin{equation*}
|\pi_1 (\bar{G}_0 (M)) | = \infty.     
\end{equation*}    
\end{thm}

\begin{proof}
We follow the proof of Theorem \ref{thm:4.1}, using $\hat K^\eta.$ 
To show (i) and (ii) in Theorem \ref{thm:5.1},
it is enough to use $\phi^* \alpha =\alpha$. 
To show (iii), we take $H=1$.  
\end{proof}

 \bibliographystyle{amsplain}
\bibliography{Paper}

\end{document}